\documentclass[12pt]{amsart}
\usepackage{amsmath, amsthm, amsfonts, stmaryrd, fullpage}

\def\FF{\mathbb{F}}

\def\QQ{\mathbb{Q}}
\def\ZZ{\mathbb{Z}}

\newcommand{\hot}{\mbox{(higher order terms)}}

\renewcommand{\setminus}{-}
\newcommand{\conts}{{\operatorname{conts}}}
\DeclareMathOperator{\chara}{char}
\DeclareMathOperator{\Tr}{Tr}

\newcommand{\F}{{\mathbb F}}

\newcommand{\Q}{{\mathbb Q}}

\newcommand{\Z}{{\mathbb Z}}

\newcommand{\Kbar}{{\overline{K}}}

\newcommand{\calL}{{\mathcal L}}

\DeclareMathOperator{\id}{id}

\DeclareMathOperator{\Aut}{Aut}

\newcommand{\GL}{\operatorname{GL}}

\newcommand{\isom}{\simeq}

\newcommand{\Intersection}{\bigcap}

\newcommand{\Union}{\bigcup}

\newcommand{\directsum}{\oplus}

\newtheorem{theorem}{Theorem}[section]
\newtheorem{cor}[theorem]{Corollary}
\newtheorem{prop}[theorem]{Proposition}
\newtheorem{lemma}[theorem]{Lemma}
\newtheorem{conj}[theorem]{Conjecture}

\theoremstyle{remark}
\newtheorem{remark}[theorem]{Remark}

\theoremstyle{definition}
\newtheorem*{hypothesis}{Hypothesis A}

\title{Orbits of automorphism groups of fields}
\author{Kiran S. Kedlaya}
\address{Department of Mathematics, Massachusetts Institute of Technology,
Cambridge, MA 02139-4307, USA}
\email{kedlaya@math.mit.edu}
\urladdr{http://math.mit.edu/\~{}kedlaya}
\author{Bjorn Poonen}
\address{Department of Mathematics, University of California, 
        Berkeley, CA 94720-3840, USA}
\email{poonen@math.berkeley.edu}
\urladdr{http://math.berkeley.edu/\~{}poonen}
\date{June 28, 2005}

\begin{document}

\begin{abstract}
We address several specific aspects of the following general question:
can a 
field $K$ have so many automorphisms that the
action of the automorphism group on the elements of $K$ has relatively few
orbits? We prove that any field which has
only finitely many orbits under its automorphism group is finite. We extend
the techniques of that proof to approach a broader conjecture, which asks 
whether
the automorphism group of one field over a subfield can have only finitely many
orbits on the complement of the subfield. Finally, we apply similar methods
to analyze the field of Mal'cev-Neumann ``generalized power series'' over a
base field; these form near-counterexamples to our conjecture when the base
field has characteristic zero, but often fall surprisingly far short in positive
characteristic.
\end{abstract}

\maketitle

Can an infinite field $K$ have so many automorphisms that the
action of the automorphism group on the elements of $K$ has only
finitely many orbits?
In Section~1, we prove that the answer is ``no'' (Theorem~\ref{thm1}),
even though the corresponding answer for division rings is probably ``yes'' 
(see Remark~\ref{remark:noncommutative}).
Our proof constructs a ``trace map'' from the given field to a
finite field, and exploits the peculiar combination of additive and
multiplicative properties of this map.

Section~2 attempts to prove a relative version of Theorem~\ref{thm1}, 
by considering, for a nontrivial extension of fields $k \subset K$,
the action of $\Aut(K/k)$ on $K$.
In this situation each element of $k$ forms an orbit,
so we study only the orbits of $\Aut(K/k)$ on $K \setminus k$.
Our Conjecture~\ref{conjecture-sec2}
asserts that if $\Aut(K/k)$ acts on $K \setminus k$
with finitely many orbits,
then $k$ and $K$ are either both finite or both algebraically closed. 
This conjecture contains Theorem~\ref{thm1} as a special case,
as one sees by taking $k$ to be the minimal subfield of $K$.
Using variants of the techniques of Section~1
(including a ``norm map'' serving as a multiplicative analogue of our
earlier ``trace map''), we prove some weaker versions of our
conjecture. For instance, under the
hypothesis of Conjecture~\ref{conjecture-sec2} and the assumption that
$k$ and $K$ are not finite, 
$k$ satisfies
Kaplansky's ``Hypothesis~A'' (Proposition~\ref{P:additive}),
and both $k$ and $K$ are radically closed
(Corollary~\ref{cor:radically closed}).

Whereas the results of Sections 1 and~2 restrict the possibilities
for fields with many automorphisms (few orbits), 
Section~3 investigates some specific candidates for fields 
that {\em could} have many automorphisms.
Specifically, we study the Mal'cev-Neumann
fields of ``generalized power series'' over a base field $k$.
If $k$ has characteristic zero or satisfies Hypothesis~A, 
then the Mal'cev-Neumann field over $k$ has relatively
few orbits under its automorphism group
(Theorem~\ref{T:orbits under Hypothesis A}), though not so few as 
to contradict Conjecture~\ref{conjecture-sec2}.
In contrast, if $k$ does not satisfy Hypothesis A, then
the Mal'cev-Neumann field over $k$ has only automorphisms given by
rescaling the series parameter (Theorem~\ref{thm2}).
The techniques used here are similar to those used in the previous
sections; indeed, the historical order of things is that we considered
the Mal'cev-Neumann fields as a source of potential counterexamples
to Conjecture~\ref{conjecture-sec2}, and the
ideas used in the proof of Theorem~\ref{thm2} led ultimately to
the proof of Theorem~\ref{thm1}.

\section{Fields with finitely many orbits}
\label{S:fields}

This section is devoted to the proof of the following theorem.
\begin{theorem} \label{thm1}
Let $K$ be a field on which the number of orbits of $\Aut(K)$
is finite. Then $K$ is finite.
\end{theorem}

\begin{remark} \label{remark:noncommutative}
The noncommutative analogue is probably false.
A related result, Corollary~2 on page~117 of \cite{Cohn},
states that given a division ring $K$ and a field $k$ contained in
its center, one can embed $K$ in a division ring $L$
such that two elements of $L$ are conjugate if and only
if they are both transcendental over $k$
or both algebraic with the same minimal polynomial over $k$.
If we take $k=\F_2$ and $K=\F_2(t)$ where $t$ is an indeterminate,
it seems as if the construction there would produce a division ring $L$
in which all elements of $L-\{0,1\}$ are conjugate;
indeed this is claimed at the bottom of page~117 of \cite{Cohn}.
But, as George Bergman pointed out to us,
it appears that this relies upon the plausible but unproved assumption
that when one forms a ``free product'' of division rings without
elements algebraic over the prime field (outside the prime field),
the resulting division ring itself has no such algebraic elements.
\end{remark}

For the rest of this section, we assume that $K$ is a field such that
\begin{center}
\framebox{The number of orbits of $\Aut(K)$ on $K$ is finite.}
\end{center}
The integral closure $k$ of the prime subfield in $K$ must be
a finite field $\FF_q$, or else $k$ alone would
contribute infinitely many orbits. 
Write $q=p^e$ where $p$ is prime and $e \in \Z_{\ge 1}$.
The $\F_p$-vector space $K$ can be made
a module over the polynomial ring $\FF_p[F]$
by setting $F \cdot \alpha = \alpha^p$ for all $\alpha \in K$.

\begin{lemma}[Bergman and Kearnes] \label{lem:trace}
There exists a unique map $\Tr\colon K \to \FF_q$
such that for any $x \in K$ and any nonzero
$P(F) \in \FF_p[F]$, there exists $y \in K$
satisfying $P(F)(y) = x - \Tr(x)$.
Moreover, $\Tr(F(x)) = F(\Tr(x))$, and $\Tr(\Tr(x)) = \Tr(x)$,
and $\Tr$ is $\F_q$-linear.
If $x$ and $y$ are in the same orbit of $\Aut(K/\FF_q)$, 
then $\Tr(x) = \Tr(y)$.
\end{lemma}
In particular (this being the case we will need), 
if $\Tr(x) = 0$, then for all $n$, there exists $y_n \in K$ such that
$y_n^{p^n} - y_n = x$. 
(This consequence of Lemma~\ref{lem:trace} could also be proved directly.)

\begin{remark}
A map like $\Tr$ exists in some other contexts; for one example, see
Section~\ref{sec:malneu}.
\end{remark}

\begin{remark}
One can make an analogous multiplicative construction: 
for every $x \in K^*$, 
there is a unique $c \in \FF_q^*$ such that $x/c$ 
has an $n$-th root for each positive integer $n$. 
We will not need this yet,
but it will come up in Section~\ref{sec:rel}.
\end{remark}

Using the interpretation of $K$ as an $\FF_p[F]$-module, 
we may deduce Lemma~\ref{lem:trace}
from the following lemma suggested by Hendrik Lenstra,
which may be of independent interest.
Let $R$ be an integral domain.
A submodule $N$ of an $R$-module $M$ is said to be \emph{characteristic} 
if $\alpha(N)=N$ for every module automorphism $\alpha$ of $M$.
For $a \in R$, define $M_a:=\{m \in M:am=0\}$.
The {\em torsion submodule} $M_t$ of $M$ is $\Union_{a \in R-\{0\}} M_a$.
The {\em divisible submodule} $M_d$ of $M$ is the set of $m \in M$
such that for every $a \in R-\{0\}$ there exists $x \in M$ with $ax=m$.

\begin{lemma}[Lenstra] \label{lem:lenstra}
Let $R$ be an integral domain, 
and let $M$ be an $R$-module with finitely many characteristic submodules.
Then $M = M_t \directsum M_d$.
Moreover there exists $c \in R - \{0\}$ annihilating $M_t$,
and for any such $c$ we have $cM=M_d$.
\end{lemma}

\begin{proof}
The submodule $M_a$ is characteristic for any $a \in R$, so there
are only finitely many possibilities for it. 
If $M_{a_1}, \dots, M_{a_n}$ is an exhaustive list of the $M_a$ 
for $a \ne 0$, put $c = a_1\cdots a_n$. 
Then $M_c = M_t$.

{}From now on, let $c$ be any nonzero element of $R$ annihilating $M_t$.
Since $R$ is a domain, $cM$ is torsion-free.
For any nonzero $a \in R$, the chain
\[
cM \supseteq acM \supseteq a^2cM \supseteq \cdots
\]
must be eventually constant, since each term is characteristic.
Choose $r$ such that $a^r cM = a^{r+1} cM$; 
then $cM = acM$ since $a^r \ne 0$ and $cM$ is torsion-free.
This holds for all nonzero $a$, so $cM \subseteq M_d$.
On the other hand, $M_d \subseteq cM_d \subseteq cM$, so $M_d=cM$.
In the exact sequence
\[
        0 \to M_c \to M \stackrel{c}\to cM \to 0,
\]
we have $M_c=M_t$, and the submodule $M_d \subseteq M$
maps isomorphically to $cM=M_d$,
so $M = M_t \directsum M_d$.
\end{proof}

\begin{proof}[Proof of Lemma~\ref{lem:trace}]
Set $M = K$ and $R = \FF_p[F]$.
Any field automorphism of $K$ is an $R$-module automorphism, 
because Frobenius commutes with all field automorphisms.
Characteristic submodules are unions of orbits under $\Aut(K)$,
so there are at most finitely many.
Moreover, $M_t = \FF_q$ because $\FF_q$ is integrally closed in $K$.
The condition characterizing $\Tr$ in Lemma~\ref{lem:trace}
says that $\id_K-\Tr$ maps $K$ into $M_d$.
If $\Tr\colon K \to \F_q$ is any map satisfying this condition,
then $\id_K$ decomposes as the sum of
$\Tr$ (which maps $K$ into $\F_q=M_t$) and $\id_K-\Tr$
(which maps $K$ into $M_d$);
thus $\Tr$ can only be the projection $M \to M_t$
in the decomposition $M \isom M_t \directsum M_d$
of Lemma~\ref{lem:lenstra}.

We now observe that the map $\Tr$ defined this way
has the claimed properties.
By construction, $\Tr$ satisfies the condition involving $P(F)$.
By Lemma~\ref{lem:lenstra}, $M_d=(F^e-1)M$, 
which is an $\F_q$-subspace of $K$.
Since $\Tr$ is a projection for a decomposition $M_t \directsum M_d$
into $\F_q$-subspaces,
it is an $\F_q$-linear map satisfying $\Tr(\Tr(x)) = \Tr(x)$.
Since $F$ maps $M_t$ and $M_d$ into themselves, $\Tr(F(x))=F(\Tr(x))$ holds.
By uniqueness, $\Tr$ is equivariant for field automorphisms,
so $\Tr$ is constant on orbits of $\Aut(K/\FF_q)$. 
\end{proof}

In the notation of the previous proof,
$F$ maps $M_t=\FF_q$ onto itself,
and maps $M_d$ onto itself,
so $K$ is perfect.
For $x \in K$ and $n \in \ZZ$, 
define $s_n(x) = \Tr(x^{1+p^n}) = \Tr(xF^n(x))$.

\begin{lemma} \label{lem:periodic}
There exists $m \in \Z_{\ge 1}$ such that 
$s_{m+n}(x) = s_n(x)$ for all $x \in K$ and $n \in \ZZ$.
\end{lemma}

\begin{proof}
Since $\Aut(K/\F_q)$ has finite index in $\Aut(K)$,
the set $S$ of $\Aut(K/\F_q)$-orbits is finite.
Let $\sim$ be the equivalence relation on $K$ induced by this partition.
For each positive integer $i$ divisible by $e$,
the map $x \mapsto F^i(x)-x$ induces a map $S \to S$.
Since $S$ is finite, there exist $i<j$ for which these maps coincide. 
(Thanks to Bergman for pointing this out, thus supplanting a
more complicated construction.)
For $n \in \Z$ and $x \in K$, we have 
\[
F^n(F^i(x) - x) \cdot (F^i(x) - x) \sim 
F^n(F^j(x) - x) \cdot (F^j(x) - x),
\]
which expands to
\begin{align*}
& F^{n+i}(x) F^i(x) - F^n(x) F^i(x) - F^{n+i}(x) x + F^n(x) x 
\\
&\sim
F^{n+j}(x) F^j(x) - F^n(x) F^j(x) - F^{n+j}(x) x + F^n(x) x.
\end{align*}
Applying $\Tr$, and using the fact that $F^i$ and $F^j$ act trivially
on the image of $\Tr$, we get
\begin{align*}
s_n(x) - s_{n-i}(x) - s_{n+i}(x) + s_n(x) 
        &= s_n(x) - s_{n-j}(x) - s_{n+j}(x) + s_n(x) \\
s_{n+j}(x) - s_{n+i}(x) - s_{n-i}(x) + s_{n-j}(x) 
        &= 0.
\end{align*}
For fixed $x$, this linear recurrence implies
that the sequence $(s_n(x))_{n \in \ZZ}$ is periodic
(since $s_n(x) \in \FF_q$ for all $n$).
The coefficients of the recurrence are independent of $x$,
so only finitely many sequences are possible,
so one can find a uniform period that works for all $x$.
\end{proof}

\begin{lemma} \label{lem:last}
Suppose $x \in K$ satisfies $\Tr(x) = 0$.
Then $\Tr(x^2) = 0$. 
If moreover $p=2$, then $\Tr(x^3) = 0$ also.
\end{lemma}

\begin{proof}
Choose $m$ as in Lemma~\ref{lem:periodic}.
Replace $m$ by a multiple if necessary, to assume that $e$ divides $m$.
Since $\Tr(x) = 0$, there exists $y \in K$ such that $F^m y - y = x$.
Now
\begin{align*}
\Tr(x^{2}) &= \Tr((y^{p^m} - y)^{2}) = \Tr(y^{2p^m}) - 2\Tr(y^{p^m+1}) + \Tr(y^2) \\
&= s_0(y)^{p^m} - 2 s_m(y) + s_0(y) = 2 s_0(y) - 2 s_m(y) = 0.
\end{align*}
In case $p=2$, we also have
\begin{align*}
\Tr(x^3) &= \Tr((y^{2^m} + y)^3) 
= \Tr(y^{3\cdot 2^m}) + \Tr(y^{2^{m+1}+1}) + \Tr(y^{2^m+2}) + \Tr(y^3) \\
&= s_1(y)^{2^m} + s_{m+1}(y) + s_{m-1}(y^2) + s_1(y)\\
&= s_{1}(y) + s_{-1}(y^2) = \Tr(y^3) + \Tr(y^3) = 0.
\end{align*}
\end{proof}

\begin{proof}[Proof of Theorem~\ref{thm1}]
We claim that $K=\FF_q$.
If not, then there exists $x \in K^*$ with $\Tr(x) = 0$. 
Let $c = \Tr(x^{-1})$.
Lemma~\ref{lem:last} implies
\begin{align*}
0 &= \Tr((x+x^{-1}-c)^2) \\
&= \Tr(x^2) + 2 \Tr(1) - 2c \Tr(x) + \Tr((x^{-1}-c)^2) = 0+2-0+0,
\end{align*}
which is a contradiction if $p \ne 2$. 
If $p=2$, applying Lemma~\ref{lem:last} repeatedly yields
\begin{align*}
0 &= \Tr((x^2+x^{-1}-c)^3) = \Tr((x^4+(x^{-2}+c^2))(x^2+(x^{-1}+c))) \\
&= \Tr(x^6) + (\Tr(x^3) + c \Tr(x^4)) + (1 + c^2 \Tr(x^2)) +
\Tr((x^{-1}-c)^3) \\
&= 0 + 0 + 0 + 1 + 0 + 0,
\end{align*}
again a contradiction.
\end{proof}

\begin{remark}
Lenstra points out an alternate argument for $p \neq 2$: 
Lemma~\ref{lem:last} 
and the identity $xy = \frac{(x+y)^2-x^2-y^2}{2}$
imply that $\ker(\Tr)$ is closed under multiplication.
Then 
\[
        K \cdot \ker(\Tr) = (\F_q + \ker(\Tr)) \ker(\Tr) \subseteq \ker(\Tr)
\]
so $\ker(\Tr)$ is an ideal of $K$.
Since $\Tr|_{\F_q}$ is not identically zero,
$\ker(\Tr)$ can be only the zero ideal, 
so $K=\F_q$.
\end{remark}

\begin{remark}
\label{R:definable subsets}
Here are two model-theoretic statements for a field $K$
that are equivalent to Theorem~\ref{thm1}:
\begin{enumerate}
\item
\label{I:definable version}
If the set of definable subsets of $K$
(i.e., $\emptyset$-definable subsets in the language of rings)
is finite, then $K$ is finite.
\item
\label{I:type version}
If the set of complete $1$-types realized by elements of $K$
is finite, then $K$ is finite.
\end{enumerate}
Statements \ref{I:definable version}
and~\ref{I:type version} are equivalent to each other,
because their hypotheses are equivalent:
the set of elements of $K$ having a given type
is by definition an intersection of definable subsets.
Statement~\ref{I:definable version} implies Theorem~\ref{thm1},
since each definable subset of $K$ is a union of $\Aut(K)$-orbits.

Finally, let us prove that Theorem~\ref{thm1} implies 
Statement~\ref{I:type version}.
Let $K$ be a field in which only finitely many complete $1$-types are realized.
By Theorem~9.14 of \cite{poizat} (with the comments preceding Theorem~9.13 of \cite{poizat}), or by Exercise~10.2.5 of~\cite{hodges},
there exists an elementary extension $L$ of $K$ that 
is strongly $\omega$-homogeneous.
Since $L$ is an elementary extension, 
every element of $L$ is of one of the finitely many types realized
by elements of $K$.
But $L$ is strongly $\omega$-homogeneous, 
so any two elements of $L$ of the same type
are in the same orbit of $\Aut(L)$.
Thus $L$ has finitely many orbits.
Applying Theorem~\ref{thm1} to $L$, we find that $L$ is finite.
So its subfield $K$ also is finite.
\end{remark}

\begin{remark}
One may ask to what extent Theorem~\ref{thm1} may be generalized to
larger classes of rings; that is, one may ask for which classes of infinite 
rings $R$ does $\Aut(R)$ always act with infinitely many orbits. 
For example, we do not know whether there exists
an infinite integral domain $R$ such that $\Aut(R)$ has 
finitely many orbits on $R$:
the proof of Theorem~\ref{thm1} seems inadequate to treat this case,
since the formal properties of $\Tr$ are satisfied, for instance,
by the constant coefficient map $\F_q[[t]] \to \F_q$.
\end{remark}

\begin{remark}
We can construct infinite commutative rings with finitely many orbits.
If $V$ is a nonzero vector space over $\FF_p$, then the ring
$R = \FF_p \directsum V$ in which $vw = 0$ for all $v,w \in V$ 
carries an action of the group $\GL(V)$,
and there are $2p$ orbits,
namely $\{a\}$ and $\{a\} + (V - \{0\})$ for all $a \in \F_p$.
\end{remark}

\begin{remark}
Here is an example of an infinite {\em reduced} (but disconnected)
commutative ring whose automorphism group 
acts with finitely many orbits. Let $C$ be the Cantor set,
and let $R$ be the ring of continuous functions from $C$ with its usual 
topology to $\FF_p$ with the discrete topology. Since every nonempty
open subset of $C$ is isomorphic to $C$ itself, the group $\Aut(C)$
of homeomorphisms from $C$ to $C$
acts transitively on the set of labeled 
partitions of $C$ into any fixed finite number of disjoint open subsets.
In particular, two elements of $R$ lie in the same orbit of $\Aut(R)$ if and
only if they have the same image (as functions to $\FF_p$); hence 
$\Aut(R)$ acts on $R$ with $2^p-1$ orbits.
\end{remark}

\section{Fields with relatively few orbits}
\label{sec:rel}

The term ``relatively'' in the section title refers to automorphisms
of one field relative to a subfield. 
The following conjecture includes Theorem~\ref{thm1}. 

\begin{conj} \label{conj} \label{conjecture-sec2}
Let $K/k$ be a nontrivial extension of fields.
Then the number of orbits of 
$\Aut(K/k)$ on $K \setminus k$ is finite
if and only if $k$ and $K$ 
are either both finite or both algebraically closed.
\end{conj}

The ``if'' part holds: for finite fields it is trivial,
and for algebraically closed fields
it follows from the theory of transcendence bases.

\begin{remark}
\label{R:model theory again}
In Remark~\ref{R:definable subsets},
we mentioned a model-theoretic strengthening of Theorem~\ref{thm1}
in which the hypothesis of finitely many $\Aut(K)$-orbits
was replaced by the weaker hypothesis of finitely many definable subsets.
Similarly we could ask whether for a model-theoretic strengthening
of Conjecture~\ref{conjecture-sec2}
in which the hypothesis is weakened to say only that
\[
        \{S \cap (K-k) : S\text{ is a subset of $K$ definable over $k$}\}
\]
is finite,
or equivalently that the set of types over $k$
realized by elements of $K-k$ is finite.
(The equivalence follows, since the set of elements of $K-k$ of a given type
over $k$ is an intersection of sets of the form $S \cap (K-k)$
with $S$ definable over $k$.)
These equivalent statements imply Conjecture~\ref{conjecture-sec2},
but it is not clear whether they are implied 
by Conjecture~\ref{conjecture-sec2}.
\end{remark}

For the rest of this section, we assume that
\begin{center}
\framebox{The number of orbits of $\Aut(K/k)$ on $K-k$ is finite.}\\
\framebox{The field $K$ is infinite.}
\end{center}
and hope to prove that $k$ and $K$ are algebraically closed.
We do not succeed, but we deduce a number of facts 
restricting the possibilities for $k$ and $K$.

\begin{prop}
\label{P:infinite k}
The field $k$ is infinite.
\end{prop}

\begin{proof}
If $k$ is finite, then the number of orbits of $\Aut(K/k)$ on 
$K = (K - k) \cup k$ is finite,
so the number of orbits of $\Aut(K)$ on $K$ is finite,
contradicting Theorem~\ref{thm1}.
\end{proof}

\begin{prop}
\label{P:relatively algebraically closed}
The field $k$ is integrally closed in $K$.
\end{prop}

\begin{proof}
Suppose the integral closure $\ell$ of $k$ in $K$ were not $k$.
Since $k$ is infinite, $\ell - k$ would be infinite,
so $K-k$ would contain infinitely many finite orbits,
a contradiction.
\end{proof}

\begin{prop}
\label{P:k and K are perfect}
Both $k$ and $K$ are perfect.
\end{prop}

\begin{proof}
Suppose they are of characteristic $p$.
Let $F$ be the Frobenius endomorphism of $K$, as in Section~\ref{S:fields}.
Let $S_n = F^n(K \setminus k) - F^{n+1}(K \setminus k)$.
The $S_n$ are disjoint, and each is a union of orbits of $\Aut(K/k)$,
so some $S_n$ is empty.
But $F$ defines a bijection $S_m \to S_{m+1}$ for each $m$,
so $S_0$ is empty.
In other words, $K \setminus k \subseteq K^p$.
Taking differences of elements, we obtain $K \subseteq K^p$,
so $K$ is perfect.
By Proposition~\ref{P:relatively algebraically closed},
$k$ also is perfect.
\end{proof}

For $x,y \in K \setminus k$,
we redefine $x \sim y$ to mean
that $x$ and $y$ belong to the same orbit of $\Aut(K/k)$.
The following lemma arose out of a discussion with Bergman.

\begin{lemma}\label{lem:plusone}
If $x \in K \setminus k$ and $c \in k$,
then $x \sim x + c$.
\end{lemma}

\begin{proof}
For $a \in k^*$ and $b \in k$, the linear map
$L_{a,b}\colon K \to K$ given by $x \mapsto ax+b$ 
permutes the additive cosets of $k$ in $K$. 
Let $G$ be the group formed by these maps;
then $G$ acts on the $\Aut(K/k)$-orbits in $K \setminus k$ 
(since $G$ commutes with the action of $\Aut(K/k)$). 
There is a normal subgroup $H$ of $G$ of finite index 
that acts trivially on these orbits.
Let $n=(G:H)$.
Then $L_{a^n,0} = L_{a,0}^n \in H$ for all $a \in k^*$.
Taking a commutator with $L_{1,1}$ shows that
$L_{1,a^n-1} \in H$ for any $a \in k^*$. 
By the following lemma, $L_{1,c} \in H$ for all $c \in k$.
By the 
definition of $H$, we get $x \sim x+c$ for all $x \in K-k$ and $c \in k$.
\end{proof}

\begin{lemma}
Let $k$ be an infinite perfect field.
Then for any positive integer $n$,
the additive group of $k$ is generated by elements of the form
$a^n - 1$ for $a \in k^*$.
\end{lemma}

\begin{proof}
Since $k$ is perfect,
we may assume without loss of generality that $n$ 
is not divisible by the characteristic of $k$.
Let $G$ be the additive subgroup of $k$ generated by elements of the form
$a^n - 1$ for $a \in k^*$. Since
\[
        (a^n - 1)(b^n - 1) = [(ab)^n - 1] - (a^n - 1) - (b^n - 1),
\]
$G$ is closed under multiplication.

Since $k$ is infinite, so is $G$,
and we can fix distinct nonzero $g_1,\dots,g_n \in G$.
Then for all but finitely many $x \in k$,
the element $\gamma_i:=(g_i x + 1)^n -1$ is in $G$ for $i=1,\dots, n$.
If we expand using the binomial theorem,
and view the $\gamma_i$ as the right-hand sides
of a system of linear equations in the ``variables'' 
$\binom{n}{1} x$, $\binom{n}{2}x^2$, \dots, $\binom{n}{n} x^n$,
then Cramer's Rule gives a formula $D \binom{n}{1} x = D'$,
where $D = \det(g_i^j)_{1 \le i,j \le n} \ne 0$,
and $D'$ is given by some polynomial in the $g_i$ and $\gamma_i$ 
with integer coefficients and no constant term.
By the previous paragraph, $D' \in G$.
Thus $Dnx \in G$ for all but finitely many $x \in k$.
Since $Dn$ is nonzero and independent of $x$,
the elements $Dnx$ exhaust all but finitely many elements of $k$.
Thus $k \setminus G$ is finite.

On the other hand, $k \setminus G$ is a union of cosets of 
the infinite group $G$,
so $k \setminus G = \emptyset$.
Hence $G=k$.
\end{proof}

\begin{lemma}
\label{L:module2}
Let $R$ be a ring, and let $M$ be an $R$-module.
For each $r \in R$, let $M_r$ be the submodule of $M$
annihilated by $r$.
Suppose that $f \in R$ is such that $M_f$ has no nonzero proper submodules.
Also suppose that there is a proper submodule $N$ of $M$
such that the sequence $(f^m(M\setminus N))_{m \ge 1}$
has only finitely many distinct sets.
Then for some $n \ge 0$, $f^n M = f^{n+1} M$
and $M$ is the direct sum of $M_{f^n}$ and $f^n M$.
\end{lemma}

\begin{proof}
The descending sequence of sets $(f^m(M \setminus N))_{m \ge 1}$ 
must stabilize,
so $f^m(M\setminus N)=f^{m+1}(M\setminus N)$ for some $m \ge 1$.
Taking the submodule generated by both sides yields $f^m M = f^{m+1} M$, 
since the submodule generated by $M\setminus N$ equals $M$.

Let $n$ be the smallest nonnegative integer such that $f^n M=f^{n+1}M$.
If $n=0$, we are done, so assume $n>0$.
Applying $f$ yields $f^{n+1}M=f^{n+2}M$ and so on,
so $f^n M=f^{n'} M$ for all $n' \ge n$.
Thus, in the exact sequence
\[
        0 \to M_{f^n} \to M \stackrel{f^n}\to f^nM \to 0,
\]
the submodule $f^nM$ of $M$ in the middle 
surjects onto $f^{2n}M=f^nM$ on the right.
In particular, $M_{f^n} + f^nM = M$.

It remains to show that $M_{f^n} \cap f^nM = 0$.
If not, let $f^e$ be the smallest power of $f$ annihilating 
$M_{f^n} \cap f^nM$;
then $f^{e-1}(M_{f^n} \cap f^nM)$ 
is a nonzero submodule of $M_f \cap f^nM$.
The hypothesis on $M_f$ implies $M_f \cap f^nM = M_f$,
so $M_f \subseteq f^nM$.

Applying $f^{n-1}$ to $M_{f^n} + f^n M = M$ yields
\[
        f^{n-1}M \subseteq M_f + f^{2n-1}M = M_f + f^n M = f^n M
\]
contradicting the minimality of $n$.
\end{proof}

\begin{remark}
\label{R:module2}
If under the hypotheses of Lemma~\ref{L:module2}
one has also $f^n M = M$,
then $M \stackrel{f}\to M$ is surjective.
(This is because $M = f^n M = f^{n+1} M = f(f^n M) = f(M)$.)
\end{remark}

\begin{prop}\label{prop:notsurj}
For each $l \ge 1$, the $l$-th power maps on $k$ and $K$ are surjective.
\end{prop}

\begin{proof}
Since $k$ is integrally closed in $K$,
it suffices to prove the result for $K$.
We may reduce to the case that $l$ is prime.
By Proposition~\ref{P:k and K are perfect},
we may assume $l$ is not the characteristic of $k$.

The hypotheses of Lemma~\ref{L:module2} 
with $R=\ZZ$, $M=K^*$, $N=k^*$, and $f=l$ (the $l$-th power map)
hold 
since each set $f^m(K^* \setminus k^*)$ 
is a union of $\Aut(K/k)$-orbits in $K \setminus k$.
Let $N\colon K^* \to (K^*)_{l^n}$ be the projection $M \to M_{f^n}$
given by the direct sum decomposition in Lemma~\ref{L:module2}.
The construction of $N$ is invariant under $\Aut(K/k)$,
so $x \sim y$ implies $N(x) \sim N(y)$, which in turn implies
$N(x) = N(y)$ because $N(x)$ and $N(y)$ are in $(K^*)_{l^n} = (k^*)_{l^n} \subseteq
k$.

By Lemma~\ref{lem:plusone},
for any $x \in K \setminus k$, we have $x^{-1} \sim x^{-1} + 1$, 
so $N(x^{-1}) = N(x^{-1} + 1)$. 
Multiplying by $N(x)$, we get $1 = N(1 + x)$. 
In other words, $N(y) = 1$ for all $y \in K \setminus k$, 
and hence for all $y \in K^*$. 
Thus $M=f^n M$.
By Remark~\ref{R:module2}, $f$ is surjective;
that is, the $l$-power map on $K^*$ is surjective.
\end{proof}

Our next proposition is an additive analogue 
of Proposition~\ref{prop:notsurj}.
Call a polynomial $P(x)$ {\em additive} 
if $P(x+y)=P(x)+P(y)$ as polynomials.

\begin{prop} \label{P:additive}
Every nonzero additive polynomial over $k$ induces surjective maps
on $k$ and~$K$.
\end{prop}

In particular, $k$ satisfies (the field-theoretic component of)
Kaplansky's ``Hypothesis A''; see the next section.

\begin{proof}
We may assume $\chara(k) = p>0$.
It suffices to consider additive polynomials $P$ of degree $>1$ 
that cannot be written as the composite of two other
additive polynomials of degree $>1$.
We will apply 
Lemma~\ref{L:module2} with $R=\FF_p[P]$ (the subring generated by $P$
in the endomorphism ring of the additive group of $K$), 
$M=K$, $N=k$, and $f=P$.
As in the proof of Proposition~\ref{prop:notsurj},
each set $f^m(M\setminus N)$ 
is a union of $\Aut(K/k)$-orbits in $K \setminus k$.
We need also to check that the kernel of $P\colon K \to K$
has no nonzero proper submodules.
This holds, because by Proposition~1.8.2 of~\cite{goss}
such a submodule $Z$ would give rise to
a nontrivial factorization $P = Q \circ R$ of additive polynomials over $k$
where $R(x)=\prod_{\alpha \in Z} (x-\alpha)$ is in $k[x]$
(each $\alpha$ is in $k$, since $k$ is integrally closed in $K$).

Let $\Tr\colon K \to \ker(P^n)$ be the projection $M \to M_{f^n}$
given by the direct sum decomposition in Lemma~\ref{L:module2}.
Again $x \sim y$ implies $\Tr(x) = \Tr(y)$.
For $x \in K \setminus k$
and $c \in \ker(P^n) \subseteq k$
we have $x \sim x+c$ by Lemma~\ref{lem:plusone}.
Applying $\Tr$ yields $\Tr(x) = \Tr(x) + c$,
so $\ker(P^n)=0$.
By Remark~\ref{R:module2}, $P$ is surjective on $K$.
Since $k$ is integrally closed in $K$, $P$ is surjective on $k$ also.
\end{proof}

\begin{prop}
\label{P:no abelian extensions}
The field $K$ has no nontrivial abelian extensions.
The same is true of~$k$.
\end{prop}

\begin{proof}
Because $k$ is integrally closed in $K$,
it suffices to prove 
that $K$ has no abelian extensions of degree $n \ge 2$.
We prove this by strong induction on~$n$.

Suppose that $n \ge 2$, and the result is known for every $n'<n$.

\bigskip
\noindent{\em Case 1: $n$ is not prime.}
The result for $n$ follows from the result for 
the prime factors of $n$.

\bigskip
\noindent{\em Case 2: $n=\chara(k)$.}
By Proposition~\ref{P:additive}, the map $x \mapsto x^n-x$ on $K$
is surjective,
so by Artin-Schreier theory, $K$ has no abelian extension of degree~$n$.

\bigskip
\noindent{\em Case 3: $n$ is a prime other than $\chara(k)$.}
Adjoining all $n$-th roots of unity to $K$ gives an abelian extension
of degree at most $\phi(n) < n$;
by the inductive hypothesis this extension is trivial.
Thus the $n$-th roots of unity are already in $K$.
By Kummer theory, all abelian extensions of $K$ of degree $n$
are contained in the field $K_n$ obtained by adjoining 
the $n$-th roots of all elements of $K$.
Proposition~\ref{prop:notsurj} implies that $K_n=K$,
so abelian extensions of $K$ of degree $n$ do not exist.
\end{proof}

\begin{cor} \label{cor:radically closed}
The field $K$ is radically closed
(that is, if $x \in \Kbar$ and $x^n \in K$ for some $n \ge 1$,
then $x \in K$).
The same is true of~$k$.
\end{cor}

\begin{cor}
If $\chara(k) = p>0$, then $k$ contains an algebraic closure of~$\FF_p$.
\end{cor}

\section{Automorphisms of Mal'cev-Neumann fields}
\label{sec:malneu}

For $k$ a field and $G$ an ordered abelian group,
the \emph{Mal'cev-Neumann field} $k((t^G))$ 
is the set of formal sums $\sum_{i \in G} c_i t^i$ 
whose {\em support} $\{ i: c_i \ne 0\}$ 
is a well-ordered subset of $G$;
multiplication is given by formal series convolution
\[
\left(\sum_{i \in G} c_i t^i \right)
\left(\sum_{i \in G} d_i t^i \right)
= \sum_{i \in G} \left( \sum_{j \in G} c_j d_{i-j} \right) t^i.
\]
This construction actually dates back to Hahn~\cite{hahn},
but the names of Mal'cev and Neumann are often associated to this field
because they generalized the construction to the case of 
a division ring $k$ and a nonabelian ordered group $G$,
in which case $k((t^G))$ is a division ring.
The elements of $k((t^G))$ are sometimes also called 
``generalized power series''.

There is a natural inclusion of fields $k \hookrightarrow k((t^G))$ 
mapping $c$ to $c t^0$.
Suppose $x = \sum_{i \in G} c_i t^i \in k((t^G))^*$.
The smallest $j$ such that $c_j \ne 0$ 
is called the {\em valuation} $v(x)$ of $x$.
For that $j$, we call $c_j t^j$ the {\em leading term} of $x$,
and call $c_j$ the {\em leading coefficient} of $x$.
Call $x$ {\em monic} if its leading coefficient is $1$.
The map $v\colon k((t^G))^* \to G$ is a valuation in the usual sense.
Define $v(0):=\infty$.
The disjoint union $G \cup \{\infty\}$ is ordered so that $g < \infty$
for all $g \in G$.

{}From now on, we take $G=\Q$.
Then $k((t^\Q))$ has an absolute value defined by $|x|:=e^{-v(x)}$
for nonzero $x$.
Let $\Aut_\conts(k((t^\Q))/k)$ be the group
of continuous automorphisms of $k((t^\Q))$ 
whose restriction to $k$ is the identity.
A continuous automorphism $\phi$ need not preserve the valuation,
but it is easy to show that for each $\phi$
there exists $r \in \Q_{>0}$ such that $v(\phi(x))=r v(x)$.

\subsection{Automorphisms in the presence of Hypothesis A}

\begin{theorem}
\label{T:substitution}
Suppose that $k$ is a field of characteristic~$0$.
For any monic $x \in k((t^\Q))^*$ of positive valuation
there exists $\phi_x \in \Aut_\conts(k((t^\Q))/k)$ mapping $t$ to $x$,
defined by ``substitution''.
\end{theorem}

\begin{proof}
We will define $\phi_x(\sum c_i t^i)$ as $\sum c_i x^i$,
but we need to make sense of the latter.

Write $x=t^m(1+\epsilon)$ where $m \in \Q_{>0}$ and $v(\epsilon)>0$.
Define $x^i=t^{mi} \sum_{n=0}^\infty \binom{i}{n} \epsilon^n$;
since $v(\epsilon^n) \to \infty$,
the series converges to an element of $k((t^\Q))$.
Next, if one substitutes this definition of $x^i$ into 
$\sum c_i x^i$, one obtains a double series of monomials in $t$
such that there are only finitely many monomials having
a given exponent, and the set of all occurring exponents
is well-ordered; this follows from the following standard lemmas.
(Here 
$S_1+\cdots+S_n:=\{\,s_1+\cdots+s_n : s_i \in S_i \text{ for all $i$}\,\}$
and $nS:= S + \cdots + S$.)
\begin{enumerate}
\item[(i)] 
If $S_1, \dots, S_n$ are well-ordered subsets of $\QQ$, then 
$S_1 + \cdots + S_n$ is well-ordered (\cite[Lemma~13.2.9(ii)]{passman}
in the key case $n=2$).
\item[(ii)]
If $S_1, \cdots, S_n$ are well-ordered subsets of $\QQ$, then
for any $x \in \QQ$, the number of $n$-tuples $(s_1, \dots, s_n) \in
S_1 \times \cdots \times S_n$ such that $s_1 + \cdots + s_n = x$ is finite
(\cite[Lemma~13.2.9(i)]{passman} in the key case $n=2$).
\item[(iii)]
If $S$ is a well-ordered subset of $\QQ \cap (0, +\infty)$, then
$\tilde{S} = \cup_{n=1}^\infty nS$
also is well-ordered; moreover, $\cap_{n=1}^\infty n\tilde{S} = \emptyset$
\cite[Lemma~13.2.10]{passman}.
\end{enumerate}
Collecting terms with the same exponent,
we obtain an element of $k((t^\Q))$, 
and we define $\phi_x(\sum c_i t^i)$ to be this element.

A similar argument shows that $\phi_x$
respects addition and multiplication.
It also acts as the identity on $k$.
Looking at leading terms shows that if $y \in k((t^\Q))^*$,
then 
\begin{equation}
\label{E:valuation}
        v(\phi_x(y))=v(x)v(y).
\end{equation}
In particular, $\phi_x$ is injective and continuous.
Also by~\eqref{E:valuation},
$k((t^\Q))$ is an immediate extension of $\phi_x(k((t^\Q)))$,
but the latter is abstractly isomorphic to $k((t^\Q))$
and hence is maximally complete (see~\cite{kap} for definitions).
Thus this immediate extension is trivial.
Hence $\phi_x$ is an automorphism.
\end{proof}

\begin{remark}
\label{R:substitution in char p} 
The proof of Theorem~\ref{T:substitution} does not work 
in characteristic $p>0$, as we now explain.
The binomial theorem does not apply to $(1+\epsilon)^i$ if $p$
divides the denominator of $i$.
Instead one must write $i = p^b q$ where $b \in \Z$
and $q \in \Q$ has denominator not divisible by $p$,
and define
\[
        (1+\epsilon)^i = \left((1+\epsilon)^q \right)^{p^b}
\]
where $(1+\epsilon)^q$ is defined using the binomial theorem,
and the map $z \mapsto z^{p^b}$ is defined termwise.
But now if $x=t-t^2$ and $y=t^{-1/p} + t^{-1/p^2} + \dots$,
then $\phi_x(y)$ makes no sense,
since a short calculation shows that
the double series that should represent it
has infinitely many terms of valuation $0$.
\end{remark}

The phenomenon in Remark~\ref{R:substitution in char p} 
was observed already by Kaplansky
in the course of his study of immediate maximal extensions 
of valued fields \cite{kap}; 
this study hinges on a key definition, which we now recall.

\begin{hypothesis}
If $k$ is a field of characteristic $p>0$
and $G$ is an ordered abelian group,
say that the pair $(k,G)$ \emph{satisfies Hypothesis~A}
is satisfied if the following two conditions hold:
\begin{enumerate}
\item
Every nonzero additive polynomial over $k$ 
induces a surjective map from $k$ to itself;
i.e., for any $a_0, \dots, a_n \in k$ not all zero and any $b \in k$,
the equation
\[
a_n x^{p^n} + \cdots + a_1 x^p + a_0 x = b
\]
has a solution $x \in k$. (In particular, $k$ is perfect.)
\item
The group $G$ is $p$-divisible, i.e., $pG = G$.
\end{enumerate}
If $G$ is omitted, we say that $k$ satisfies Hypothesis~A if the first
condition above holds. 
As discussed in \cite[pp.~20--21]{kap2},
Whaples \cite{whaples} proved that
$k$ satisfies Hypothesis~A 
if and only if $k$ has no finite extension of degree divisible by $p$.
If instead $k$ has characteristic~$0$,
then by convention, $k$ and $(k,G)$ satisfy Hypothesis~A.
\end{hypothesis}

We now have the following generalization of Theorem~\ref{T:substitution}.

\begin{theorem}
\label{T:automorphisms under Hypothesis A}
Suppose that $k$ is a field satisfying Hypothesis~A.
For any monic $x \in k((t^\Q))^*$ of positive valuation
there exists $\phi_x \in \Aut_\conts(k((t^\Q))/k)$ mapping $t$ to $x$.
\end{theorem}

\begin{proof}
Let $k(t^\Q)$ be the subfield of $k((t^\Q))$ generated by $k$ and $t^i$
for all $i \in \Q$.
For $i \in \Q$, define $x^i \in k((t^\Q))$ as in the proof
of Theorem~\ref{T:substitution}, using the modifications
outlined in Remark~\ref{R:substitution in char p}.
Let $k(x^\Q)$ be the subfield of $k((t^\Q))$ generated by $k$ and $x^i$
for all $i \in \Q$.
If we forget the embeddings into $k((t^\Q))$,
then there is a $k$-isomorphism $k(t) \to k(x)$
mapping $t$ to $x$;
this extends to a $k$-isomorphism $k(t^\Q) \to k(x^\Q)$
mapping $t^i$ to $x^i$ for each $i \in \Q$.
Now $k((t^\Q))$ is a maximally complete immediate extension 
of both $k(t^\Q)$ and $k(x^\Q)$ (see~\cite{kap} for definitions),
so by \cite[Theorem~5]{kap}, the $k$-isomorphism $k(t^\Q) \to k(x^\Q)$
extends to a continuous automorphism $k((t^\Q)) \to k((t^\Q))$
(still mapping $t$ to $x$).
\end{proof}

\begin{theorem}
\label{T:orbits under Hypothesis A}
Suppose that $k$ is a radically closed field satisfying 
Hypothesis~A.
Define
\begin{align*}
        S_0 &= \{\,y \in k((t^\QQ))^* : v(y)>0\,\} \\
        S_\infty &= \{\,y \in k((t^\QQ))^* : v(y)<0\,\} \\
        S_c &= c + S_0 \qquad \text{for $c \in k$.}
\end{align*}
Then the $S_c$ for $c \in k$ and $S_\infty$
are all the orbits of $\Aut_\conts(k((t^\QQ))/k)$ 
on $k((t^\QQ)) \setminus k$.
\end{theorem}

\begin{proof}
For any group homomorphism $\lambda\colon \Q \to k^*$,
define $\psi_\lambda \in \Aut_\conts(k((t^\Q))/k)$ by
\[
        \psi_\lambda \left(\sum c_i t^i \right) = \sum \lambda(i) c_i t^i.
\]
Given $i>0$ and $a \in k^*$,
we can find $\lambda$ such that $\psi_\lambda(t^i)=a t^i$,
since $k$ is radically closed.
Thus every element of $S_0$
is in the same orbit as a monic element of $S_0$.
By Theorem~\ref{T:automorphisms under Hypothesis A},
every monic element of $S_0$ is in the same orbit as $t$.
Thus $S_0$ is contained in an orbit.

On the other hand, $S_0$ is preserved by each continuous automorphism 
of $k((t^\Q))$, 
since
\[
        S_0 = \{\,x \in k((t^\Q))^*: x^n \to 0 \text{ as } n \to \infty\,\}.
\]
Thus $S_0$ is an orbit.

The maps $x \mapsto x+c$ for $c \in k$ and $x \mapsto x^{-1}$
are $\Aut_\conts(k((t^\Q))/k)$-equivariant bijections from
$k((t^\Q)) \setminus k$ to itself,
so they map orbits to orbits.
Thus each $S_c$ is an orbit, and $S_\infty$ is an orbit.
Their union is all of $k((t^\Q)) \setminus k$,
so they are all the orbits.
\end{proof}

\begin{remark}
If we used $\Aut(k((t^\Q))/k)$ in place of $\Aut_\conts(k((t^\Q))/k)$,
the orbits could be even larger.
For example, if $k$ is algebraically closed,
then $k((t^\Q))$ is algebraically closed,
so $k((t^\Q)) \setminus k$ consists of one orbit under $\Aut(k((t^\Q))/k)$.
\end{remark}

\subsection{Automorphisms in the absence of Hypothesis A}

Now, in the spirit of \cite[Section~5]{kap}, 
we consider what happens when Hypothesis~A fails in the field aspect;
we find that $k((t^{\QQ}))$ has very few endomorphisms over $k$.

\begin{theorem} \label{thm2}
Let $k$ be a perfect field not satisfying Hypothesis~A.
Then the endomorphisms of $k((t^\Q))$ over $k$ are the maps
of the form
\[
\sum_i c_i t^i \mapsto 
\sum_i \lambda(i) c_i t^{ri}
\]
where $\lambda\colon \QQ \to k^*$
is a group homomorphism and $r \in \QQ_{>0}$.
\end{theorem}

In particular all endomorphisms of $k((t^\Q))$ are automorphisms,
and they are all continuous.

\begin{cor}
If $k$ is finite, then the endomorphisms of $k((t^\Q))$ over $k$ 
are the maps of the form 
$\sum c_i t^i \mapsto \sum c_i t^{ri}$
where $r \in \QQ_{>0}$.
\end{cor}

\begin{proof}
If $q=\#k$, then $x^q-x=1$ has no solution in $k$,
so $k$ does not satisfy Hypothesis~A.
Apply Theorem~\ref{thm2} and observe that
every group homomorphism $\lambda\colon \QQ \to k^*$ is trivial.
\end{proof}

We will deduce Theorem~\ref{thm2} from a slightly more general result,
Theorem~\ref{thm3} below.
Let $\Tr\colon k((t^\QQ)) \to k$ be the ``trace'' map carrying a series
$x = \sum c_i t^i$ to its constant coefficient $c_0$. 
For any field $K$ with $k \subseteq K \subseteq k((t^\Q))$,
let $K^{\Tr} = \{\, x \in K : \Tr(x) = 0 \,\}$;
this is a $k$-subspace of $K$.
Let $p$ be the characteristic of $k$.

\begin{lemma}
\label{L:bijective P}
If $k$ is perfect, then each nonzero additive polynomial $P$ over $k$
maps $k((t^\QQ))^{\Tr}$ bijectively to itself.
\end{lemma}

\begin{proof}
The additive polynomials $x \mapsto x^p$ and $x \mapsto ax$ for $a \in k$
map $k((t^\QQ))^{\Tr}$ into itself.
Any additive polynomial can be built from these using composition
and addition, so $P$ maps $k((t^\QQ))^{\Tr}$ into itself.
Since $k$ is integrally closed in $k((t^\QQ))$,
each $P$ acts injectively on $k((t^\QQ))^{\Tr}$.

It remains to show that 
$P \colon k((t^\QQ))^{\Tr} \to k((t^\QQ))^{\Tr}$
is surjective.
The result is true for $x \mapsto x^p$ 
so we may reduce to the case in which $P$ is separable.
By additivity, it suffices to solve $P(x)=b$
in the following two cases.

\bigskip
\noindent{\em Case 1: $b$ has only positive exponents.}

Then $v(b)>0$.
Since the lowest degree monomial in $P$ has degree~$1$,
there exists a formal power series solution
\[
        c_0 b + c_1 b^p + c_2 b^{p^2} + \cdots 
\]
with coefficients in $k$; 
this converges to an actual solution to $P(x)=b$.

\bigskip
\noindent{\em Case 2: $b$ has only negative exponents.}

Since $k$ is perfect, one can solve for coefficients $c_i \in k$ making
\begin{equation}
\label{E:formal solution}
        c_n b^{1/p^n} + c_{n+1} b^{1/p^{n+1}} + \cdots
\end{equation}
a formal solution, where $\deg P = p^n$.
Since $b$ has only negative exponents,
the same is true for each $b^{1/p^m}$.
Moreover, given $\epsilon>0$,
only finitely many of the $b^{1/p^m}$
contribute monomials with exponents more negative
than $-\epsilon$.
Thus the series~\eqref{E:formal solution} 
makes sense as an element of $k((t^\Q))^{\Tr}$;
moreover, it represents a solution to $P(x)=b$.
\end{proof}

\begin{theorem} \label{thm3}
Let $k$ be a perfect field not satisfying Hypothesis~A.
Suppose $K$ is a field such that $k(t^\Q) \subseteq K \subseteq k((t^\QQ))$ 
and $P(K^{\Tr}) = K^{\Tr}$ for each nonzero additive polynomial $P$
over $k$.
Then the $k$-homomorphisms $s\colon K \to k((t^\Q))$ 
are the maps of the form
\[
\sum_i c_i t^i \mapsto 
\sum_i \lambda(i) c_i t^{ri}
\]
where $\lambda\colon \QQ \to k^*$
is a group homomorphism and $r \in \QQ_{>0}$.
\end{theorem}

\begin{remark}
In the special case where $k$ is a finite field $\FF_q$,
a slight modification of our proof (left to the reader) 
shows that the hypothesis $P(K^{\Tr}) = K^{\Tr}$ 
need be assumed only for $P(x)=x^q-x$.
\end{remark}

\begin{remark}
Theorem~\ref{thm3} applies, for instance, when 
$k$ is a perfect field not satisfying Hypothesis~A
and $K$ is the integral closure of $k(t)$ or $k((t))$ in $k((t^\Q))$.
(Both of these integral closures can be described fairly explicitly: 
see \cite{me,me2}.)
\end{remark}

The rest of this section will be devoted to proving Theorem~\ref{thm3}.
We thus assume for the remainder of this section that
\begin{center}
\framebox{The field $k$ is perfect and does not satisfy Hypothesis~A.} \\
\framebox{The field $K$ satisfies 
                $k(t^\Q) \subseteq K \subseteq k((t^\QQ))$.}\\
\framebox{For each nonzero additive polynomial $P$ over $k$, 
                we have $P(K^{\Tr}) = K^{\Tr}$.}\\
\framebox{We have a $k$-homomorphism $s\colon K \to k((t^\Q))$.}
\end{center}

We first need some auxiliary results in the spirit of the proof of
Theorem~\ref{thm1}.

\begin{lemma} \label{L:intersection}
We have $\Intersection_{P} P(K)=K^{\Tr}$,
where the intersection is taken over all nonzero additive polynomials $P$
over $k$.
\end{lemma}

\begin{proof}
Let $I$ be the intersection.
Each $P(K)$ is a subgroup of $K$, 
and multiplication by an element of $k$
permutes these subgroups, 
so $I$ is a $k$-subspace of $K$.
Since $P(K) \supseteq P(K^{\Tr}) = K^{\Tr}$ for each $P$,
we have $I \supseteq K^{\Tr}$.
But $K^{\Tr}$ has codimension~$1$ in $K$,
and $I \ne K$ because $k$ does not satisfy Hypothesis~A.
Thus $I=K^{\Tr}$.
\end{proof}

\begin{lemma}\label{L:defining trace}
We have $\Tr(s(x))=\Tr(x)$ for all $x \in K$.
\end{lemma}

\begin{proof}
Since $s$ acts trivially on $k$, 
it suffices to consider the case $x \in K^{\Tr}$.
Then $x \in \Intersection_P P(K)$ by Lemma~\ref{L:intersection},
so $s(x) \in \Intersection P(s(K)) \subseteq \Intersection P(k((t^\Q)))$,
and the latter equals $k((t^\Q))^{\Tr}$,
by Lemma~\ref{L:intersection} applied to $k((t^\Q))$.
\end{proof}

\begin{lemma}\label{L:positive valuation}
Suppose $x \in K^*$ and $\Tr(x)=0$.
Then $v(x)>0$ if and only if $\Tr\left(\frac{x^p}{x^p-x}\right)=0$.
\end{lemma}

\begin{proof}
Since $\Tr(x)=0$, we have $v(x)\ne 0$.
If $v(x)>0$, then $v\left(\frac{x^p}{x^p-x}\right) = (p-1)v(x)>0$,
so $\Tr\left(\frac{x^p}{x^p-x}\right)=0$.
If $v(x)<0$, then $x^p$ and $x^p-x$ have the same leading term,
so $\Tr\left(\frac{x^p}{x^p-x}\right)=1$.
\end{proof}

Lemmas \ref{L:defining trace} and~\ref{L:positive valuation}
imply the following:

\begin{cor} \label{C:sign}
If $x \in K^*$, then $v(x)$ and $v(s(x))$ have the same sign.
\end{cor}

The map $\calL\colon k((t^\Q))^* \to k^* t^\Q$ 
that returns the leading term of a series is a group homomorphism,
so the map $\Q \mapsto k^* t^\Q$ defined by $i \mapsto \calL(s(t^i))$
must have the form $i \mapsto \lambda(i) t^{ri}$
for some homomorphism $\lambda \colon \Q \to k^*$
and some $r \in \Q$.
Corollary~\ref{C:sign} shows that $r>0$.
By composing $s$ with an automorphism of $k((t^\Q))$ 
of the type described in Theorem~\ref{thm2},
we reduce to the following case:
\begin{center}
\framebox{For all $i \in \Q$, the leading term of $s(t^i)$ is $t^i$.}
\end{center}
We now hope to prove that $s(x)=x$ for all $x \in K$.

\begin{lemma} \label{L:power}
Under the boxed assumptions, we have $s(t)=t$.
\end{lemma}

\begin{proof}
If not, then for some $b \in \Q_{>0}$ and $c \in k^*$, we have
\[
        s(t) = t(1 + ct^b + \hot).
\]
Write $b=p^e b'$ where $e$ is the $p$-adic valuation of $b$.
Choose a large negative integer $\ell$ not divisible by $p$,
and set $j=b'/\ell$.
Thus $j<0$, the $p$-adic valuation of $j$ is $0$,
and $b/j = p^e \ell \in \Z[1/p]$.
By choosing $|\ell|$ large enough, we may assume also that $0<j+b$.

We compute $s(t^j)$ by raising $s(t)$ to an integer power,
and then taking an integer root.
Since the latter integer is prime to $p$,
and since $s(t^j)$ has leading coefficient $1$ by hypothesis,
we obtain
\begin{align*}
        s(t^j) &= t^j(1  + jc t^b + \hot) \\
        &= t^j  + jc t^{j+b} + \hot).
\end{align*}

Let $n$ be a positive integer greater than $-e$.
For any $x \in k((t^\Q))$, define $h(x)$ to be the $y \in k((t^\Q))$
such that $y^{p^n}+y=x-\Tr(x)$.
It is unique by Lemma~\ref{L:bijective P},
which also describes how to compute it.
Moreover, if $x \in K$, then $h(x) \in K$, by hypothesis.
Lemma~\ref{L:defining trace} implies that $s(h(x))=h(s(x))$.
We compute
\begin{align*}
        h(t^j) &= t^{j/p^n} + 
                \text{ (other terms with smaller negative exponent) } \\
        s(h(t^j)) &= h(s(t^j)) \\
                &= t^{j/p^n} + \text{ (terms with negative exponent) } 
                        + jct^{j+b} + \hot \\
                &= t^{-a} + \text{ (terms with negative exponent) } 
                        + jct^{pma} + \hot,
\end{align*}
where $a=-j/p^n \in \Q_{>0}$ and $m:=-(j+b)p^{n-1}/j \in \Z_{>0}$.
In the multinomial expansion for $s(h(t^j))^{1+pm}$,
any product involving at least one of the terms of $s(h(t^j))$ 
with positive exponent
will have exponent at least $pm(-a) + 1(pma) = 0$.
Moreover, there is exactly one product in the multinomial expansion
with exponent exactly $0$, namely
\[
        \binom{1+pm}{1} \left(t^{-a}\right)^{pm} (jct^{pma})^1 = (1+pm)jc,
\]
which is nonzero in $k$.
Thus $\Tr\left(s(h(t^j))^{1+pm}\right) \ne 0$.

On the other hand, $h(t^j)^{1+pm}$ has only terms with negative exponent,
so $\Tr\left(h(t^j)^{1+pm}\right)=0$.
This contradicts Lemma~\ref{L:defining trace}.
\end{proof}

\begin{cor}\label{C:powers}
For every $i \in \Q$, we have $s(t^i)=t^i$.
\end{cor}

\begin{proof}
This follows from $s(t)=t$ and the assumption that 
the leading coefficient of $s(t^i)$ is~$1$.
\end{proof}

Now, for any $x \in K$ and $l \in \Q$,
\begin{align*}
        \Tr \left( t^{-l} s(x) \right)
        &= \Tr \left( s(t^{-l} x) \right)
                &&\text{(by Corollary~\ref{C:powers})} \\
        &= \Tr \left(t^{-l} x \right) 
                &&\text{(by Lemma~\ref{L:defining trace}).}
\end{align*}
In other words, the coefficient of $t^l$ in $s(x)$ 
equals the coefficient of $t^l$ in $x$.
This holds for all $l$, so $s(x)=x$.
This completes the proof of Theorem~\ref{thm3}.

\subsection*{Acknowledgments}
This paper is the result of a discussion with George Bergman,
Keith Kearnes, Hendrik Lenstra, and Thomas Scanlon;
their contributions are noted throughout the text.
We thank also Anand Pillay and Alex Wilkie
for discussions concerning
Remarks \ref{R:definable subsets} and~\ref{R:model theory again}.
Kedlaya was partially supported by a National Science Foundation
postdoctoral fellowship (grant DMS-0071597) and NSF grant DMS-0400727.
Poonen was partially supported by NSF grants DMS-9801104 and DMS-0301280
and a Packard Fellowship; 
he thanks also the Isaac Newton Institute for hosting a visit in
the summer of 2005.


\begin{thebibliography}{Ke2}

\bibitem[C]{Cohn}
P. M. Cohn,
\textit{Skew field constructions},
London Math.\ Soc.\ Lecture Note Series {\bf 27},
Cambridge Univ.\ Press, Cambridge, 1977.

\bibitem[G]{goss}
D. Goss, 
\textit{Basic Structures of Function Field Arithmetic},
Ergebnisse der Mathematik und ihrer Grenzgebiete {\bf 35}, 
Springer-Verlag, 1996.

\bibitem[Ha]{hahn}
H. Hahn, \"Uber die nichtarchimedische Gr\"o\ss ensysteme (1907),
reprinted in \textit{Gesammelte Abhandlungen}~I, Springer-Verlag, 1995.

\bibitem[Ho]{hodges}
W. Hodges, \textit{Model Theory}, Cambridge Univ.\ Press, Cambridge, 1993.

\bibitem[Ka1]{kap}
I. Kaplansky, Maximal fields with valuations, \textit{Duke Math.
J.} \textbf{9} (1942), 303--321.

\bibitem[Ka2]{kap2}
I. Kaplansky, \textit{Selected Papers and Other Writings}, Springer-Verlag,
1995.

\bibitem[Ke1]{me}
K.S. Kedlaya, The algebraic closure of the power series field in positive
characteristic, \textit{Proc. Amer. Math. Soc.}
\textbf{129} (2001), 3461--3470.

\bibitem[Ke2]{me2}
K.S. Kedlaya, Finite automata and algebraic extensions of function fields,
preprint, \texttt{arXiv: math.AC/0410375}; to appear in
\textit{J.\ Th\'eor.\ Nombres Bordeaux}.

\bibitem[Pa]{passman}
D.S. Passman, \textit{The Algebraic Structure of Group Rings},
Wiley, 1977.

\bibitem[Po]{poizat}
B. Poizat, \textit{A Course in Model Theory} (translated by M. Klein),
Springer-Verlag, New York, 2000.

\bibitem[W]{whaples}
G. Whaples, Galois cohomology of additive polynomials and $n$-th power
mappings of fields, \textit{Duke Math. J.} \textbf{24} (1957), 143--150.

\end{thebibliography}
\end{document}